\documentclass{amsart} 
\usepackage{amssymb}
\usepackage{longtable}
\usepackage{tikz}
\usetikzlibrary{calc,positioning, shapes.geometric}

\usepackage{enumitem}
\usepackage{url}
\usepackage[colorinlistoftodos,prependcaption,textsize=tiny]{todonotes}

\theoremstyle{definition}

\theoremstyle{remark}

\numberwithin{equation}{section}

\begin{document}

\title[Primitive Complete Normal Basis Generators]{Computational Results on the Existence of Primitive Complete Normal Basis Generators}

\author[D. Hachenberger]{Dirk Hachenberger}
\address{Institut f\"ur Mathematik der Universit\"at Augsburg, D-86135 Augsburg}
\curraddr{}
\email{hachenberger@math.uni-augsburg.de}
\thanks{}

\author[S. Hackenberg]{Stefan Hackenberg}
\address{Angerstr. 30, D-86179 Augsburg}
\curraddr{}
\email{mail@stefan-hackenberg.de}
\thanks{}

\subjclass[2010]{Primary 11T30, 12E20, Secondary 05B25, 51D20}

\date{16 December 2019}

\dedicatory{}

%======================================================================================================= 

%% own-macros m
\let\cal\mathcal

\newcommand{\cfq}{\overline{\mathbb{F}}_q}
\newcommand{\F}{\mathbb{F}} 
\newcommand{\K}{\mathbb{K}} 
\newcommand{\fq}{\mathbb{F}_q} 
\newcommand{\fp}{\mathbb{F}_p} 
\newcommand{\fqm}{\mathbb{F}_{q^m}} 
\newcommand{\fqd}{\mathbb{F}_{q^d}} 
\newcommand{\fqe}{\mathbb{F}_{q^e}} 

\newcommand{\fqn}{\mathbb{F}_{q^n}} 
\newcommand{\fqk}{\mathbb{F}_{q^k}} 
\newcommand{\fql}{\mathbb{F}_{q^\ell}}
\newcommand{\N}{\mathbb{N}}
\newcommand{\Z}{\mathbb{Z}} 
\newcommand{\C}{\mathbb{C}} 
\newcommand{\ft}[1]{\mathbb{F}_{q^{#1}}} 
\newcommand{\tr}{{\rm Tr}} 

\newcommand{\G}{\mathcal G}

\newcommand{\rad}{{\rm rad}} 
\newcommand{\eps}{\varepsilon}
\newcommand{\opp}{{\rm \footnotesize opp}}

\newcommand{\Ord}{{\rm Ord}} 
\newcommand{\Tr}{{\rm Tr}} 
\newcommand{\No}{{\rm No}} 
\newcommand{\Nm}{{\rm Nm}} 

\newcommand{\ord}{{\rm ord}}

\newcommand{\grp}{{\rm grp}}
\newcommand{\id}{{\rm id}}

\newcommand{\sord}{{\rm subord}}
\renewcommand{\char}{{\rm char}}
\renewcommand{\span}{{\rm span}}
\renewcommand{\max}{{\rm max}}
\renewcommand{\min}{{\rm min}}
\newcommand{\equi}{\Leftrightarrow}
\newcommand{\impli}{\Rightarrow}
\newcommand{\lcm}{{\rm lcm}}
\renewcommand{\ker}{{\rm kern}}
\newcommand{\imag}{{\rm im}}
\newcommand{\supp}{{\rm supp}}
\newcommand{\sgn}{\mathop{\rm sgn}}
\renewcommand{\log}{{\rm log}}
\renewcommand{\det}{{\rm det}} 
\renewcommand{\deg}{{\rm deg}}
\renewcommand{\dim}{\mathop{\rm dim}}
\newcommand{\DIV}{\, {\tt div}\, }
\newcommand{\MOD}{\, {\tt mod}\, }

%======================================================================================================= 

\begin{abstract} 
We present computational results which strongly support a conjecture of Morgan and Mullen (1996), which states that for every extension $E/F$ 
of Galois fields there exists a primitive element of $E$ 
which is completely normal over $F$. 
\end{abstract}

\maketitle

%---------------------------------------------------------------------
\section{Introduction} \label{introduction} 
%---------------------------------------------------------------------
\noindent 
To every prime power $q>1$ and every integer $n\geq 1$ there corresponds (up to isomorphism) a unique 
extension $E/F$ of Galois fields: the ground field $F=\fq$ is the finite field with $q$ elements and 
$E=\fqn$ is its $n$-dimensional extension, a field with $q^n$ elements. 
It is an important open problem, whether for every pair $(q,n)$ there exists a 
\textbf{primitive} element in the corresponding field extension $E/F$, which also satisfies the property of being 
\textbf{completely normal} over $F$: 
\begin{itemize} 
\item a \textit{primitive} element of $E$ is a generator of the (cyclic) multiplicative group of $E$;  
\item 
an element $w\in E$ is called \textit{normal} over $F$, if its conjugates under the (cyclic) Galois group of $E/F$, that is, $w$, $w^q$, ..., $w^{q^{n-1}}$,  
constitute an $F$-basis of $E$; 
 \item if $w \in E$ simultaneously is 
normal over $K$ for \textit{every} intermediate field $K$ of $E/F$, then $w$ is called \textit{completely normal} over $F$. 
(Every positive divisor $d$ of $n$ gives rise to a 
unique intermediate field of $E/F$ of the form $\fqd$, and vice versa.)  
\end{itemize} 

\noindent 
Of course, if $w\in E$ is a primitive completely normal element over $F$, then so are its conjugates. 

In 1996, it has been conjectured by Morgan and Mullen \cite{Mo-Mu-1996} that for every pair $(q,n)$ there does 
exist a primitive completely normal element (for short: a PCN-element) in the corresponding field extension. 
A proof of this conjecture would generalize two fundamental theorems: 
\begin{itemize} 
\item 
first, the \textit{Primitive Normal Basis Theorem} (Lenstra and Schoof, 1987, \cite{Le-Sc-1987}), which states that 
for every extension $E/F$ of Galois fields there exists a primitive element of $E$ which is normal over $F$ (a PN-element for short); 
\item on the other hand,   
the \textit{Complete Normal Basis Theorem} (Blessenohl and Johnsen, 1986, \cite{Bl-Jo-1986}), which says that every extension of 
Galois fields admits a completely normal element (a CN-element for short). 
\end{itemize} 

\noindent 
Morgan and Mullen \cite{Mo-Mu-1996} based their conjecture 
on a computational search, whose range comprises all pairs $(p,n)$, where $p< 100$ is a prime number and where $p^n< 10^{50}$. 
They have also computed the exact number of 
all completely normal and all primitive completely normal elements for the pairs $(q,n)$ listed in Table \ref{tab:enumeration-IM-GM}, where $q$ is a prime power.

\begin{table}\label{tab:enumeration-IM-GM} 
\caption{Complete enumeration of CN- and of PCN-elements: the range of Morgan and Mullen \cite{Mo-Mu-1996}.}
\begin{tabular}{|r|l||r|l||r|l| } 
\hline 
$q$ &   & $q$ &  &$q$&  \\ \hline 
$2$ & $2 \leq n \leq 18$ & $7$ & $2 \leq n \leq 6$  &  
$3$ & $2 \leq n \leq 12$ \\ $8$ & $2 \leq n \leq 5$ & 
$4$ & $2 \leq n \leq 9$  & $9$ & $2 \leq n \leq 5$ \\  
$5$ & $2 \leq n \leq 8$  & & & & \\ 
\hline 
\end{tabular} 
\end{table}

\vspace{1ex} 
\noindent 
The aim of the present work is to use the available structural results on CN-elements, from Hachenberger \cite{Ha-1997}, 
together with a skillful implementation, based on Hackenberg \cite{Ha-2015c}, in order to extend the computational results of Morgan and Mullen 
enormously. Our main contributions are as follows. 

\vspace{1ex} 
\noindent 
{\bf Computational Result 1.} 
Let $\mathcal G$ denote the set of all integers $n\geq 1$ such that for 
\textit{every} prime power $q>1$ there exists a PCN-element in the corresponding extension $\fqn/\fq$. 
Then:  
$$n \in {\mathcal G} \ {\mbox{ for every $n$ with $1\leq n \leq 202$.}}$$

\vspace{1ex} 
\noindent 
{\bf Computational Result 2.} A monic polynomial $f(x)\in \fq[x]$ is called a \textbf{PCN-polynomial}, 
if it is irreducible over $\fq$ and if its roots are primitive and completely normal elements for $\fqn$ over $\fq$, where 
$n=\deg(f)$. For every prime number $p<10\,000$ and for every degree $n$ such that $p^n < 10^{80}$ we have determined a PCN-polynomial 
of degree $n$ over the prime field $\fp$.

\vspace{1ex} 
\noindent 
{\bf Computational Result 3.} The exact number of all CN-elements and of all PCN-elements 
for $\fqn$ over $\fq$ are determined for the pairs $(q,n)$ listed in Table \ref{tab:enumeration-DH-SH-1} and Table \ref{tab:enumeration-DH-SH-2}.  

\begin{table}
\caption{Complete enumeration of CN- and of PCN-elements: extended ranges (1).}
\label{tab:enumeration-DH-SH-1} 
\begin{tabular}{|r|l||r|l||r|l|} 
\hline 
$q$ &                  & $q$  & & $q$ &  \\ \hline 
2 & $2 \leq n \leq 31$ & 
3 & $2 \leq n \leq 20$ & 
4 & $2 \leq n \leq 14$ \\ 
5 & $2 \leq n \leq 12$ & 
7 & $2 \leq n \leq 11$ & 
8 & $2 \leq n \leq 9$ \\ 
9 & $2 \leq n \leq 9$  & 
11 & $2 \leq n \leq 7$  & 
13 & $2 \leq n \leq 7$  \\ 
16 & $2 \leq n \leq 7$  & 
17 & $2 \leq n \leq 7$  & 
19 & $2 \leq n \leq 7$  \\ 
23 & $2 \leq n \leq 7$  & 
25 & $2 \leq n \leq 6$  & 
27 & $2 \leq n \leq 4$  \\ 
29 & $2 \leq n \leq 6$  & 
32 & $2 \leq n \leq 4$  & 
37 & $2 \leq n \leq 6$  \\ 
41 & $2 \leq n \leq 6$  & 
43 & $2 \leq n \leq 6$  & 
121 & $2 \leq n \leq 4$  \\ 
169 & $2 \leq n \leq 4$  & 
361 & $2 \leq n \leq 3$  & 
529 & $2 \leq n \leq 3$  \\ 
841 & $2 \leq n \leq 3$  & 
961 & $2 \leq n \leq 3$  & 
1369 & $n=2$  \\ 
1681 & $n=2$  & 
1369 & $n=2$  & 
1849 & $n=2$  \\  
\hline 
\end{tabular} 
\end{table} 

\vspace{2ex} 
\begin{table} 
\caption{Complete enumeration of CN- and of PCN-elements: extended ranges (2).}
\label{tab:enumeration-DH-SH-2} 
\begin{tabular}{|r|l||r|l||r|l|} 
\hline 
$n$ & $q$ prime power   & $n$ & $q$ prime power  &  $n$ & $q$ prime power  \\ \hline 
3 & $2 \leq q \leq 961$ & 
4 & $2 \leq q \leq 243$ & 
6 & $2 \leq n \leq 43$ \\ 
\hline 
\end{tabular} 
\end{table}

\noindent 
For the basic theory of finite fields we refer to Lidl and Niederreiter \cite{Li-Ni-1983}, and to the forthcoming monograph 
Hachenberger and Jungnickel \cite{Ha-Ju-2020}. The latter contains proofs, both, of the primitive and of the complete normal basis theorem.

%---------------------------------------------------------------------
\section{Preliminary remarks} \label{notation} 
%---------------------------------------------------------------------
\noindent Before we are going to describe our strategies which enabled us to achieve our computational results, 
we like to comment on the present status of the Morgan-Mullen-Conjecture and on some further results concerning 
primitivity and normality of finite field elements. 

We first fix some useful notation. 
Throughout, $p$ denotes the characteristic of the underlying fields. 
Consider a positive integer $n$. 

\begin{itemize} 
\item We write $n=p^a n'$, where 
$n'$ is not divisibly by $p$, that is, $n'$ is the \textbf{$p$-free part} of $n$. 
\item The \textbf{order of $q$ modulo $n'$}, denoted by $\ord_{n'}(q)$, 
is the least integer $k\geq 1$ such that $q^k \equiv 1 \MOD n'$. 
\item Finally, $\rad(n')$ denotes the \textbf{radical} of $n'$, that is, the product over all distinct prime divisors of $n'$. 
\end{itemize} 

It will also be convenient to use the following abbreviations:  
\begin{itemize}
\item $P_n(q)$ for the number of primitive elements of $\fqn$; 
\item $N_n(q)$ for the number of normal elements of $\fqn$ over $\fq$; 
\item $PN_n(q)$ for the number of primitive normal  elements of $\fqn$ over $\fq$; 
\item $CN_n(q)$ for the number of completely normal  elements of $\fqn$ over $\fq$; 
\item $PCN_n(q)$ for the number of primitive completely normal elements of $\fqn$ over $\fq$. 
\end{itemize} 

\noindent 
For instance, when $q=2$ and $n=6$ one has 
\[P_6(2)=36, \ N_6(2)=24, \ PN_6(2)=18, \ CN_6(2)=12, \ PCN_6(2)=6.\]  

\noindent 
Of course, $P_n(q)=\varphi(q^n-1)$, where $\varphi$ is Euler's totient function. 
The additive ($q$-)analogon is $N_n(q)=\phi_q(x^n-1)$, where $\phi_q$ counts the units of the 
polynomial residue ring $\fq[x]/(x^n-1)$. 
In the trivial case, where $n=1$, we obviously have 
$N_1(q)=CN_1(q)$ and $P_1(q)=PN_1(q)$, hence $CN_1(q)=q-1$ and $PCN_1(q)=\varphi(q-1)$ (for every $q$).   

%--2.1-----------------------------------------------------------------------------------------------------------------
\subsection{$n$ prime} 
%--2.1-----------------------------------------------------------------------------------------------------------------
By the definition of complete normality, 
$CN_r(q)=N_r(q)=\phi_q(x^r-1)$, and therefore $PCN_r(q)=PN_r(q)$ for every prime number $r$ (independently from $q$). 

Especially for the case $n=2$ it is well known that 
$N_2(q)=CN_2(q)$ and $P_2(q)=PN_2(q)$, hence $CN_2(q)=\phi_q(x^2-1)$ and $PCN_2(q)=\varphi(q^2-1)$ (for every $q$); see \cite[Proposition 13.1.1]{Ha-Ju-2020}, for instance. 

%--2.2-----------------------------------------------------------------------------------------------------------------
\subsection{Completely basic extensions} 
%--2.2-----------------------------------------------------------------------------------------------------------------
Following a notion of Faith \cite{Fa-1957}, a pair $(q,n)$, as well as the corresponding field extension are called \textbf{completely basic}, 
if every normal element is already completely normal, that is, $CN_n(q)=N_n(q)=\phi_q(x^n-1)$ and therefore $PCN_n(q)=PN_n(q)$. 
According to \cite[Theorem 3.1]{Ha-2001a} (see also Blessenohl \cite{Bl-1990}, as well as Blessenohl and Johnsen \cite{Bl-Jo-1991}), 
the following holds: 

\begin{quote} 
\textit{Theorem} 2.2.1. The pair $(q,n)$ is completely basic, if and only if  
for every prime divisor $r$ of $n$, the number ${\rm ord} _{(n/r)'}(q)$ is not divisible by $r$.  
\end{quote} 

\noindent 
This implies that $(q,r^2)$ is completely basic for every prime $r$ and every $q$. Furthermore,  
$(q,p^m)$ is completely basic for every power of the characteristic $p$ of $\fq$.  

%--2.3-----------------------------------------------------------------------------------------------------------------
\subsection{Cubic and quartic extensions} 
%--2.3-----------------------------------------------------------------------------------------------------------------
The case $n=3$ is the first one, where the determination of $PN_n(q)$ becomes a nontrivial task. 
A (general) exact formula for $PN_n(q)$ is not known for $n\geq 3$, and it is very unlikely that a \textit{simple} formula can be found (given it exists at all). 
Nevertheless, for cubic ($n=3$) and quartic ($n=4$) extensions, strong lower bounds for $PN_n(q)$ have been derived in 
Hachenberger \cite{Ha-2015, Ha-2019} by geometric considerations. In particular, when $n=4$, the following hold; 
\cite[Theorem 1.3]{Ha-2019} and \cite[Remark 1.4]{Ha-2019}: 

\begin{quote} 
{\em Theorem} 2.3.1. Under the assumption that 
$q^2+1$ is a prime number if $q$ is even, and that $\frac{1}{2}(q^2+1)$ is a prime number if $q$ is odd, one has  
\[ PN_4(q) = \left\{ 
\begin{array}{rl} 
(q-1)(q-3)\cdot \varphi(q^2-1) & {\mbox{ if $q \equiv 1\mod 4$,}} \\ 
(q-1)^2\cdot \varphi(q^2-1)  & {\mbox{ if $q \equiv 3\mod 4$,}} \\ 
q(q-1) \cdot \varphi(q^2-1) & {\mbox{ if $q \equiv 0\mod 2$.}} 
\end{array} \right. \]  
%is proved in \cite{Ha-2019}. 
\end{quote} 

\begin{quote} 
{\em Theorem} 2.3.2. Assume that $q$ is a Mersenne prime (which requires $q\equiv 3 \bmod 4$(, then 
$$PN_4(q)= (2q-2) \cdot \varphi(q-1)\cdot \varphi(q^2+1).$$ 
If $q+1$ is a Fermat prime (which requires that $q$ is even), then 
$$PN_4(q)=(q-1)\cdot \varphi(q-1)\cdot \varphi(q^2+1).$$
\end{quote}

%--2.4-----------------------------------------------------------------------------------------------------------------
\subsection{Extensions of degree $6$} 
%--2.4-----------------------------------------------------------------------------------------------------------------
Because of the above, the first degree, where the property of \textit{completeness} becomes meaningful is $n=6$. 
A study of $6$-dimensional extensions under a projective geometric point of view, providing lower bounds for $PN_6(q)$ and $PCN_6(q)$,  
is in preparation. %\cite{Ha-2016}. 

%--2.5-----------------------------------------------------------------------------------------------------------------
\subsection{Regularity} 
%--2.5-----------------------------------------------------------------------------------------------------------------
\noindent Starting with \cite{Ha-1997} (see also \cite{Ha-2013}) there have been achieved 
various results concerning the structure of completely normal elements, which led to a proof of the Morgan-Mullen-Conjecture for the special, but quite large class of \textit{regular} field extensions; \cite{Ha-2001a, Ha-2014}: 

\begin{quote} 
{\em Theorem} 2.5.1.
Assume that the pair $(q,n)$ is \textbf{regular}, which means that 
$n$ and $\ord_{\rad(n')}(q)$ are relatively prime. Then there exists a PCN-element in the corresponding extension of Galois fields. 
\end{quote} 

\noindent 
The class of regular extensions comprises (but is not restricted to) the class of all prime power 
extensions. 

\begin{quote} 
{\em Definition} 2.5.2.  
A positive integer $n$ is called \textbf{universally regular}, if $(q,n)$ is regular for every prime power $q>1$. 
\end{quote} 

\noindent 
As a consequence, with $\mathcal G$ as defined in the introducion (see Computational Result 1), we have $n\in {\mathcal G}$ whenever $n$ is universally regular. 
As remarked above, any prime power is universally regular. In fact, 
$n$ is universally regular, whenever $r$ does not divide $s-1$ for any two distinct prime divisors $r$ and $s$ of $n$. The list 
\[ \begin{array}{l} 15, 33, 35, 45, 51, 65, 69, 75, 77, 85, 87, 91, 95, 99 \\ 
                   115, 119, 123, 133, 135, 141, 143, 145, 153, 159, 161, 175, 177, 185, 187 \end{array} \] 
comprises all universally regular numbers $\leq 200$ which are not prime powers.

%--2.6-----------------------------------------------------------------------------------------------------------------
\subsection{Lower bounds for $CN_n(q)$ and $PCN_n(q)$} 
%--2.6-----------------------------------------------------------------------------------------------------------------
An exact formula for $CN_n(q)$ is known for the class of regular pairs (resp. regular extensions), see \cite{Ha-1997}. 
It is conjectured, see \cite{Ha-1997, Ha-2013}, that 
$$CN_n(q) \geq (q-1)^{n'} \cdot q^{(p^a-1)n'} $$
for \textit{all} pairs $(q,n)$, where $n=p^an'$ as above. Moreover, it is conjectured that equality holds, if and only if 
$n'$ divides $q-1$, in which case every normal element of $\fqn$ over $\fq$ already is 
completely normal in that extension. 
This bound is known to be true for all regular pairs (once more \cite{Ha-1997, Ha-2013}) 
and it is additionally supported by our computational enumerations. 

Some nontrivial 
lower bounds for $PCN_n(q)$ are provided in \cite{Ha-2010} for the case where $n$ is a prime power.
%But since it is not even known whether $PCN_n(q)>0$ for all $(q,n)$, % (this is just a reformulation of the Morgan-Mullen Conjecture), 
%an explicit formula for this function cannot be expected at all. 

%--2.7-----------------------------------------------------------------------------------------------------------------
\subsection{An asymptotic result} 
%--2.7-----------------------------------------------------------------------------------------------------------------
By \cite[Theorem 2]{Ha-2016}, for every fixed $n$, the quotient 
$PCN_n(q)/P_n(q)$ converges to $1$ as $q$ tends to infinity. 
This gives a strong asymptotical evidence for the Morgan-Mullen Conjecture.

%--2.8-----------------------------------------------------------------------------------------------------------------
\subsection{Primitive completely normal elements for large $q$}  
%--2.8-----------------------------------------------------------------------------------------------------------------
Theorem 1 of \cite{Ha-2016} provides a further sufficient condition for the existence of a PCN-element.  

\begin{quote} 
{\em Theorem} 2.8.1.  Assume that $$q \geq \frac{(t(n)-1) \cdot (\ln(2) + n \ln(q))}{\ln(2)},$$ 
where $t(n):=\sum_{d|n}d$ is the sum of all positive divisors of $n$, and where $\ln$ denotes the natural logarithm; then there exists a PCN-element in $\fqn$ over $\fq$. 
\end{quote} 
This is used to settle the asymptotic result mentioned in Subsection 2.7 and in order to show that $PCN(q,n)>0$ whenever $q\geq n^{7/2}$ and $n\geq 7$, or when $q\geq n^3$ and $n\geq 37$ (see \cite{Ha-2016}).   
Recently, the latter could be improved considerably by Garefalakis and Kapatenakis \cite{GaKa-2018}: 

\begin{quote} 
{\em Theorem} 2.8.2. One even has $PCN(q,n)>0$ whenever $q>n'$, where $n'$ (as above) is the $p$-free part of $n$. 
\end{quote}

%---------------------------------------------------------------------
\section{Sufficient Existence Criteria} 
\label{sufficient} 
%---------------------------------------------------------------------
\noindent 
The aim of the present section is to explain our strategy in order to achieve our Computational Result 1. 
It is based on a variety of sufficient
number theoretical conditions for the existence of a primitive completely normal 
element in $\fqn$ over $\fq$. 
The principal idea rests on the following simple observation:  

\begin{quote} 
assume that $U_{(q,n)}$ is an upper bound for $q^n-CN_n(q)$, that is, for the number of elements of $\fqn$ that are \textit{not} 
completely normal over $\fq$, and assume that $L_{(q,n)}$ is a lower bound for $P_n(q)$. 
If $L_{(q,n)} > U_{(q,n)}$, then clearly $PCN_n(q)>0$. 
\end{quote} 

\noindent 
In what follows, we are going to work out a concrete version of this. 
% idea, which in fact behaves fairly good on the ranges 
%we have tested. 

%--3.1-----------------------------------------------------------------------------------------------------------------
\subsection{Lower bounds for the number of primitive elements} 
%--3.1-----------------------------------------------------------------------------------------------------------------
The following elementary lower bound for the number of all primitive elements of $\fqn$ 
is considered in \cite{Ha-2016}:   
\begin{equation} \label{eqn:Pnq-1}
P_n(q) >  \frac{\ln(2)\cdot (q^n-1)}{\ln(2) + n\ln(q)} .  
\end{equation} 
The following sharper bound (used in our present approach) can be found in Rosser and Schoenfeld \cite{Ro-Sc-1962}: 

\begin{equation} \label{eqn:Pnq-2}
P_n(q) \geq \frac{q^n-1}{e^\gamma \cdot \ln(\ln(q^n-1))
      + \frac{3}{\ln(\ln(q^n-1))}}  =: L_{(q,n)}, \end{equation}  
where $e$ is Euler's number and $\gamma$ is the Euler-Mascheroni constant.

%--3.2-----------------------------------------------------------------------------------------------------------------
\subsection{A lower bound for the number of completely normal elements} 
%--3.2-----------------------------------------------------------------------------------------------------------------
In order to tighten the complete normality condition we first introduce the following 

\begin{quote} 
{\em Definition} 3.2.1.
Let $E/F$ be the extension of Galois fields corresponding to the pair $(q,n)$. 
A subset $\mathcal D$ of $\{d\in \N: d\mid n\}$ is called \textbf{$(q,n)$-essential} (or $E/F$-\textbf{essential}), provided that 
$w \in E$ is completely normal over $F$ if and only if 
$w$ is normal in $E/\fqd$ for every $d\in {\mathcal D}$.  % \qed
\end{quote}

\noindent 
Obviously, $\{d\in \mathbb{N}: d\mid n, d\not=n\}$ is $(q,n)$-essential for any $q$, since any nonzero element of $E$ is normal over $E$. 
Of course, we are interested to derive $E/F$-essential sets which are as small as possible. 
For this purpose, we are going to produce a specific non-trivial $E/F$-essential set, 
denoted by ${\mathcal D}^*$, which relies on \cite[Theorem 15.5]{Ha-1997}. 

\begin{quote} 
{\em Proposition} 3.2.2.
Let $E/K$ be an extension of Galois fields with degree $m$ and with $|K|=Q$. Let $r$ be a prime divisor of $m$, and 
let $L$ be the intermediate field of $E/K$ with degree $r$ over $K$. 
Then the following two assertions are equivalent: 
\begin{enumerate} 
\item[(i)] every normal element of $E/K$ is normal in $E/L$;  
\item[(ii)] $r$ does not divide $\ord_{(m/r)'}(Q)$.  %\qed
\end{enumerate} 
\end{quote} 

\noindent 
Next, suppose we are given a pair $(q,n)$ and a proper divisor $d$ of $n$. 
Let $r$ be a prime divisor of $n/d$ and consider $e:=dr$. 
If $\ord_{(n/e)'}(q^d) \not\equiv 0 \MOD r$, then  Proposition 3.2.2 can be applied with $K=\fqd$ and $L=\fqe$. 
This suggests to introduce the following directed graph associated with $(q,n)$. 

\begin{quote} 
{\em Definition} 3.2.3. The \textbf{CN-digraph} $\Gamma$ associated to the pair $(q,n)$ is defined as follows: \begin{itemize} 
\item the set of vertices of $\Gamma$ is the set of all positive divisors $d$ of $n$ with $d\not=n$;  
\item two distinct vertices $d$ and $e$ are connected by an arc, directed from $d$ to $e$ 
(for short: $d\rightarrow e$), provided the 
following two conditions are satisfied: 
\begin{enumerate} 
\item $d$ divides $e$ and $r:=\frac{e}{d}$ is a prime number;  
\item the order of $q^d$ modulo $(\frac{n}{e})'$ is not divisible by $r$. %\qed
\end{enumerate} 
\end{itemize} 
\end{quote} 

\noindent 
This immediately leads to the following 

\begin{quote} 
{\em Proposition} 3.2.4.  
Let $\Gamma$ be the CN-digraph associated to the pair $(q,n)$.  
Define ${\cal D}^*$ to be the set of all vertices of $\Gamma$ having indegree $0$, that means,  
$t \in {\mathcal D}^*$ if and only if there is no divisor $s$ of $n$ such that $s\rightarrow t$ is an arc of $\Gamma$. 
Then ${\mathcal D}^*$ is $(q,n)$-essential. % \qed
\end{quote} 

\noindent For example, when $q=3$ and $n=20$, one has ${\mathcal D}^*=\{1,2,4\}$. 

\vspace{1ex} 
\noindent 
Given \textit{some} $(q,n)$-essential set, we next determine a lower 
bound for the 
number of completely normal elements in the corresponding Galois field extension.

\begin{quote} 
{\em Proposition} 3.2.5.  
Let $\cal D$ be some $(q,n)$-essential subset of divisors of $n$. 
Furthermore, let 
\[ U_{(q,n)} := \sum_{d \in {\cal D}} \Big(\sum_{a|\frac{n}{d}} \mu\left(\tfrac{n}{da}\right) q^{da} \ - \ \phi_{q^d}(x^{\frac{n}{d}}-1) \Big) , \] 
where $\mu$ denotes the M\"obius function. 
Then $CN_n(q)\geq q^n-U_{(q,n)}$. 

\vspace{1ex} \noindent 
{\em Proof.}  
For $d\in {\mathcal D}$, let $G_d$ be the set of all $w\in E$ such that $\fqd(w)=E$; in 
other words, $G_d$ is the set of elements of $E$ which are not contained in a field $K$ with $\fqd \subseteq K \subseteq E$ and $K\not=E$. 
Furthermore, let $N_d$ denote the set of all $w\in E$ that are normal in $E/\fqd$.  
Since $\mathcal D$ is $E/F$-essential, 
the set $C$ of all $w\in E$ that are completely normal over $F$ is equal to 
$C=\bigcap_{d \in {\mathcal D}} N_d$. 
On the other hand, assuming that $w$ is a generator of $E/F$ (that is $w\in G_1$), which is not completely normal over 
$F$, gives that there is a $d\in {\mathcal D}$ such that $w \not\in N_d$, hence 
$w \in G_d \setminus N_d$, since $N_d \subseteq G_d$ for every $d$ (observe that $G_1 \subseteq G_d$ for all $d$). 
This implies $w \in \bigcup_{d \in {\cal D}} (G_d \setminus N_d)$, and therefore, the set $E\setminus C$ of all 
$w\in E$ which are not completely normal over $F$ satisfies 
\[ |E\setminus C| \leq \sum_{d \in {\cal D}} (|G_d| - |N_d|) =:U_{(q,n)}. \] 
Finally,  by a fundamental result of finite field theory (see for instance  \cite{Li-Ni-1983}), one has  
\[ |G_d| = \sum_{a|\frac{n}{d}} \mu\left(\tfrac{n}{da}\right) q^{da} \] 
and $|N_d| = \phi_{q^d}(x^{n/d}-1)$ (for all $d$). 
This gives the bound as claimed. \qed
\end{quote} 

\noindent 
Altogether, as explained in the introduction of this section, 
\eqref{eqn:Pnq-1} and \eqref{eqn:Pnq-2} together with Proposition 3.2.5 provide a sufficient existence criterion for 
$PCN_n(q)$ to be positive, explicitely: 

\begin{equation} \label{eqn:sufficient-1}
\frac{q^n-1}{e^\gamma \cdot \ln(\ln(q^n-1))
      + \frac{3}{\ln(\ln(q^n-1))}} > \sum_{d \in {\cal D}} \Big(\sum_{a|\frac{n}{d}} \mu\left(\tfrac{n}{da}\right) q^{da} \ - \ \phi_{q^d}(x^{\frac{n}{d}}-1) \Big). 
      \end{equation}

\vspace{1ex} 
\noindent 
We shall derive an alternative sufficient criterion, next. 
In fact, it is adopted from, and improves the correponding result in \cite{GaKa-2018}; the improvement rests on the fact that we 
work with the non-trivial $(q,n)$-essential sets which are based on the CN-graphs rather than the trivial one, $\{d\in \mathbb{N}: d\mid n, d\not=n\}$, which in fact is never optimal. 

Throughout, let $\omega=\omega(q^n-1)$ denote the number of all distinct prime divisors of $q^n-1$. 
Let $\mathcal D$ be some set which is $(q,n)$-essential. 
For every $d\in {\mathcal D}$, let 
$\Omega_d = \Omega_d(x^{n/d}-1)$ 
denote the number of distinct monic divisors of $x^{n/d}-1$ that are irreducible over $\F_{q^d}$. 
The following formula is well known: 
$$\Omega_d = \sum_{t \mid (n/d)'} \frac{\varphi(t)}{\ord_t(q^d)}.$$ 
Finally, for every $d\in {\mathcal D}$, let 
$$\Theta_d = \frac{\phi_{q^d}(x^{(n/d)'}-1)}{q^{d\cdot (n/d)'}}.$$ 
Generalizing the criterion (7), respectively (11) of \cite{GaKa-2018}, with respect to $\mathcal D$, we obtain: 

\begin{quote} 
\textit{Proposition} 3.2.6. Sufficient for $PCN_n(q)$ to be positive is 
the condition 
\begin{equation} \label{eqn:sufficient-2} 
CN(q,n) > q^{n/2} \cdot (2^\omega -1) \cdot \prod_{d \in {\mathcal D}} (\Theta_d \cdot 2^{\Omega_d}). 
\end{equation} 
Moreover, since always $\Theta_d <1$, the following is sufficient as well, where $U_{(q,n)}$ is as in Proposition 3.2.5: 
\begin{equation} \label{eqn:sufficient-3}
q^n-U_{(q,n)} \geq q^{n/2} \cdot 2^\omega \cdot 2^{\sum_{d \in {\mathcal D}} \Omega_d} 
\end{equation} 
Finally, using the bound 
$$2^\omega \leq 4514.7 \cdot q^{n/8},$$ 
established in Lemma 3.2 of \cite{GaKa-2018}, it would be sufficient to have 
\begin{equation} \label{eqn:sufficient-4}
q^n-U_{(q,n)} \geq 4514.7 \cdot  q^{5n/8} \cdot 2^{\sum_{d \in {\mathcal D}} \Omega_d},  
\end{equation} 
where, again, $U_{(q,n)}$ is taken from Proposition 3.2.5. % \qed
\end{quote}

%--3.3-----------------------------------------------------------------------------------------------------------------
\subsection{The strategy to establish Computational Result 1}  
\label{sufficient-result-1} 
%--3.3-----------------------------------------------------------------------------------------------------------------
\noindent 
Our strategy underlying the Computational Result 1 can now be summarized as follows: 

\begin{enumerate} 
\item[$\circ$] Suppose a fixed degree $n$ is given. 
Because of the result of Garefalakis and Kapatenakis \cite{GaKa-2018}, which here is Theorem 2.8.2, and because of Subsection 2.5, we only need to consider prime powers $q$ such that $q<n'$ and $(q,n)$ not regular. 
(Observe that $q=n'$ cannot happen.) 
\item[$\circ$] Given such a $q$, 
\begin{enumerate} 
\item determine first the $(q,n)$-essential set ${\mathcal D}^*$ resulting from the CN-digraph as provided in Definition 3.2.3;  
\item based on this, determine $U_{(q,n)}$ as in Proposition 3.2.5., as well as $\sum_{d \in {\mathcal D}^*} \Omega_d$. 

\vspace{1ex} 
\end{enumerate} 
\item[C1] Test, whether Inequality \eqref{eqn:sufficient-1} is satisfied. \vspace{1ex} 
\item[C2] If not, test, whether Inequality \eqref{eqn:sufficient-4} is satisfied. \vspace{1ex} 
\item[C3] If not, replace the factor $2^{\sum_{d \in {\mathcal D}^*} \Omega_d}$ in \eqref{eqn:sufficient-4} by 
$\prod_{d \in {\mathcal D}^*} (\Theta_d \cdot 2^{\Omega_d})$ and test, whether 
the condition 
$$q^n-U_{(q,n)} \geq 4514.7 \cdot  q^{5n/8} \cdot \prod_{d \in {\mathcal D}^*} (\Theta_d \cdot 2^{\Omega_d})$$ 
is satisfied. \vspace{1ex} 
\item[C4] If this still fails, determine the exact value of $\omega$; this requires the prime power factorization of 
$q^n-1$. Check now, whether Inequality \eqref{eqn:sufficient-3} is satisfied. \vspace{1ex} 

\item[C5] If this is not the case, then consider Inequality \eqref{eqn:sufficient-2} with the left hand side replaced by $q^n-U_{(q,n)}$, 
that is 
$$q^n-U_{(q,n)} > q^{n/2} \cdot (2^\omega -1) \cdot \prod_{d \in {\mathcal D}} (\Theta_d \cdot 2^{\Omega_d}).$$ 

\item[C6] If this attempt also fails, then verify the existence of a PCN-element in the current field extension by searching for a 
PCN-polynomial; a task which is explained in detail in the forthcoming section, and which of course is used to establish our Computational Result 2. 
\end{enumerate}

%---------------------------------------------------------------------
\section{Determination of PCN-polynomials}  
\label{PCN-polynomials}  
%---------------------------------------------------------------------
\noindent Recall from the discussion of our strategy in the last section that 
the concrete search for a PCN-element in $E=\fqn$ over $F=\fq$ has become necessary after the pair $(q,n)$ has failed 
all sufficient conditions provided in the last section. 
Also, in view of our Computational Result 2, we need to 
setup an explicit model for the extension field $E$ and search for a PCN-polynomial. 

%--4.1-----------------------------------------------------------------------------------------------------------------
\subsection{Modelling finite field extensions} 
%--4.1-----------------------------------------------------------------------------------------------------------------
Assume that $\fqn$ has characteristic $p$ and let $q=p^e$. Then $\fqn$ has degree $en$ over its prime field 
$\fp$ (of residues modulo $p$), and therefore $\fqn$ can be obtained as 
a residue ring $\fp[x]/(f)$, where $f(x)\in \fp[x]$ is some monic polynomial with degree 
$en$ which is irreducible over $\fp$. 
It is well known (see for instance \cite{Li-Ni-1983}) that the number of such polynomials $f$ is equal to  
\[ \frac{1}{en} \cdot \sum_{d\mid en} \mu\left(\tfrac{en}{d}\right) p^d. \] 
%An approach for finding such a polynomial $f$ which has as less non-vanishing coefficients as possible 
%is considered in Subsection 5.5. 
%%Alternatively, in the case where $e>1$, one may first describe $\fq$ as a residue ring 
%%$\fp[x]/(g)$ with $g(x) \in \fp[x]$ monic and irreducible over $\fp$ with degree $e$, and obtain $\fqn$ as 
%%residue ring $\fq[y]/(h)$ for some monic $h(y) \in \fq[y]$ which is irreducible over $\fq$ with degree $n$. 
%%We have always chosen the first way, that is, a direct representation of $E$ over its prime subfield. 

After this is done, any field element $v \in \fqn$ corresponds to a unique 
polynomial $a(x)\in \fp[x]$ with degree strictly less that $en$, namely $v=a(x) + (f)$, and the arithmetic in $\fqn$ is performed modulo $f(x)$ 
(and modulo $p$). 
Of course, $x+(f)$ is the canonical candidate to test for primitivity and complete normality, first. 
These tests are described in the forthcomming two subsections.

%--4.2-----------------------------------------------------------------------------------------------------------------
\subsection{Testing complete normality} 
%--4.2-----------------------------------------------------------------------------------------------------------------
We have to start with some preliminaries, for which we refer to \cite{Ha-1997}.   
Consider again the extension $E/F$ of Galois fields, corresponding to the pair $(q,n)$. 
The Frobenius automorphism $\sigma: E \rightarrow E$, $w \mapsto w^q$ generates the 
(cyclic) Galois group of $E/F$. Its minimal polynomial is equal to $x^n-1$. 
The \textbf{$q$-order} of $w \in E$, denoted by $\Ord_q(w)$, is the monic polynomial $g(x)\in F[x]$ of least degree 
such that 
$w$ is annihilated by the $F$-endomorphism $g(\sigma)$ (for short: $g(\sigma)w=0$).  
The $q$-order of $w$ divides $x^n-1$, and equality occurs, if and only if $w$ is normal over $F$. 
The condition that $w$ is completely normal over $F$ can therefore be phrased as 
$\Ord_{q^d}(w)=x^{n/d}-1$ for every divisor $d$ of $n$. More oeconomically, using Proposition 3.2.4, we have 

\begin{quote} 
\textit{Proposition} 4.2.1. 
An element $w\in \fqn$ is completely normal over $\fq$ if and only if 
$\Ord_{q^d}(w)=x^{n/d}-1$ for every divisor $d \in {\mathcal D}$, where $\mathcal D$ is some $(q,n)$-essential set, for instance 
the essential set ${\mathcal D}^*$ arising from the CN-digraph associated to $(q,n)$. 
\end{quote} 

\noindent 
Back to our model from Subsection 4.1, suppose we are given a concrete element $w\in \fqn=\fp[x]/(f)$, say $w=x+(f)$. 
Let ${\mathcal D}^*$ be as in Proposition 4.2.1, and let $d\in {\mathcal D}^*$. 
\begin{itemize} 
\item  We factorize the polynomial $x^{n/d}-1$ over $\fqd$, in order to get its distinct monic irreducible divisors $g_1(x), ..., g_t(x)$ 
 (over $\fqd$), and for $i=1,...,t$ let  
 $G_i(x):=(x^{n/d}-1)/g_i(x)$ be the corresponding cofactors. 
 \item Then $w$ is normal over $\fqd$, if and only if $G_i(\sigma^d)w\not=0$ for all $i$; the latter just means that 
 $w$ is not contained in any of the maximal $\sigma^d$-invariant $\fqd$-subspaces of $\fqn$. 
\end{itemize} 

\noindent If this holds for all $d\in {\mathcal D}^*$, then $w$ is a CN-element for $\fqn$ over $\fq$.  
 
 \vspace{1ex} \noindent 
 We shall mention that Morgan and Mullen \cite{Mo-Mu-1996} used a different (complete) normality test: 
 consider a divisor $d$ of $n$; then $w$ is normal in $\fqn$ over $\fqd$ if and only if 
 \[ \gcd\Big( x^{\frac{n}{d}}-1, \sum_{i=0}^{\frac{n}{d}-1} w^{q^{di}} x^{\frac{n}{d}-i} \Big) = 1. \] 
 %\todo{Ist es so, dass unser Normalit\"atstest schneller ist?} 
 In \cite{Mo-Mu-1996} this is carried out for all $d\mid n$ with $d\not=n$. But even when restricting this gcd-test to 
 divisors $d$ from  ${\mathcal D}^*$, we made the experience that the strategy for testing complete normality explained first performs faster.

%--4.3-----------------------------------------------------------------------------------------------------------------
\subsection{Testing primitivity} 
%--4.3-----------------------------------------------------------------------------------------------------------------
Suppose that the given element $w\in \fqn$ has been identified to be completely normal. 
Then the factorization of $q^n-1$ (already obtained in Step C4 of the strategy explained in the last section) 
can be used to check whether $w$ is primitive, which is done in analogy to the performance of the CN-test explained in Subsection 4.2:  
let $r_1, ...,r_k$ be all the distinct prime divisors of $q^n-1$, and let $R_i=(q^n-1)/r_i$ be their corresponding cofactors;  
then $w$ is primitive if and only if $w^{R_i}\not=1$, because then, $w$ is not contained in any of the maximal subgroups of 
the (cyclic) multiplicative group of $\fqn$. 

Of course, the square-and-multiply technique is essential when determining $w^{R_i}$. 
%, and as well by evaluating the terms 
%$G_i(\sigma^d)w$ during the CN-test. 

%--4.4-----------------------------------------------------------------------------------------------------------------
\subsection{Finding absolute PCN-polynomials} 
%--4.4-----------------------------------------------------------------------------------------------------------------
Assume next (after setting up $\fqn$ as $\fp[x]/(f)$) that the canonical candidate $x+(f)$ turned out not to be a PCN-element. 
Then, in principle, one can search through $\fqn$ until a PCN-element 
$v=a(x) + (f)$ is found by varying $a(x) \in \fp[x]$ with degree less that $en$. 
In order to identify $v$, one would then require the model parameter $f(x)$ along with the 
polynomial $a(x)$. 

In accordance with Morgan and Mullen, we have chosen the following different approach: 
instead of fixing the model and changing the polynomial $a(x)$, it is more oeconomical to change the model parameter $f(x)$ until 
the canonical candidate $w=x+(f)$ turns out to be a PCN-element, in which case only $f(x)$ has to be tabulated. 
If $w$ even is a PCN-element for the $en$-dimensional extension $\fqn$ over $\fp$, then $f(x)$ is called an \textbf{absolute PCN-polynomial}. 

\begin{quote} 
We have arranged our computations in such a way that we always determine 
\textit{absolute} PCN-polynomials. 
\end{quote} 

\noindent A further look at the tables of Morgan and Mullen \cite{Mo-Mu-1996} 
motivates the search for (absolute) PCN-polynomials with as few non-vanishing
coefficients as possible. This is respected by the following definition of 
a (total) \textbf{polynomial order} on the set of all monic polynomials of $\fp[x]$ with a fixed degree. 

\begin{quote} 
\textit{Definition} 4.4.1. 
Let first $f(x) = x^m+a_{m-1}x^{m-1}+\ldots+a_0$ be a monic polynomial with degree $m$ over the field $\fp$ of residues modulo the prime $p$. 
\begin{itemize} 
\item The \textbf{support} of $f(x)$ is the index set of its non-vanishing coefficients: 
  $\supp(f):=\{ i: a_i \not=0 \}$, where of course $a_m=1$, and $|\supp(f)|$ is the 
  \textbf{Hamming-weight} of $f(x)$. 
  \item Assume that $\supp(f)=\{i_1,...,i_k\}$ with $i_1<i_2<...<i_k$. 
  Then $I(f) := i_1|i_2|\ldots |i_k$ represents $\supp(f)$ as a word over the alphabet $\{0,1,...,m\}$ 
  in ascending order. As for the concrete coefficients of $f(x)$, we consider 
  the word $C(f):=a_{i_k}|\ldots |a_{i_2}|a_{i_1}$ over the alphabet $\fp$, which is given by 
  the canoncial residue system $\{0,1,...,p-1\}$. 
\end{itemize} 
Next, let $f(x)$ and $g(x)$ be two distinct monic polynomials with degree $m$ over $\fp$. 
  Then $f(x)$ is said to be \textbf{smaller} than $g(x)$, 
  denoted as $f\prec g$, provided one of the following conditions is satisfied: 
  \begin{enumerate}
    \item $|\supp(f)| < |\supp(g)|$;  
    \item $|\supp(f)| = |\supp(g)|$ and $I(f)$ is lexicographically smaller than $I(g)$;  
    \item $|\supp(f)| = |\supp(g)|$ and $I(f)=I(g)$ and $C(f)$ is lexicographically smaller than $C(g)$. 
  \end{enumerate}
The last comparison relies on the natural order of $\{0,1,\ldots,p-1\}$.   
\end{quote}

\noindent When searching through the set of all monic polynomials from $\fp[x]$ with degree $en$, increasing with respect to $\prec$, 
some coefficients may be restricted as follows: 
let 
$f(x) = x^{en} + a_{en-1}x^{en-1}+\ldots+a_0 \in \fp[x]$ 
be an absolute PCN-polynomial, and let $w$ be some of its roots in $\fqn$; 
  
  \begin{itemize} 
  \item then $(-1)^{e n} a_0$ is equal to the $(\fqn,\fp)$-norm of $w$ (this is the product of all its $\fp$-conjugates), and it is therefore a primitive element 
  of the prime field $\fp$; 
  \item similar, on the additive side, $a_{en-1}$ is equal to the 
  $(\fqn,\fp)$-trace of $w$ (which is the sum of all its $\fp$-conjugates of $w$), and is therefore non-zero. 
  \end{itemize} 

\noindent 
With these two restrictions in mind, 
the smallest (cf. Definition 4.4.1)  polynomials to be considered are \textbf{trinomials}:  
\[ x^{en} + \alpha x^{en-1} + \beta. \] 
In fact, we have detected plenty of absolute PCN-trinomials, and therefore offer the following conjecture. 

\begin{quote} 
\textit{Conjecture} 4.4.2. 
For every integer $m\geq 2$ there is a bound $T_m$ with the following property:  
for every prime $p\geq T_m$ there exists an absolute PCN-trinomial of degree $m$ over $\fp$. 
\end{quote}

%---------------------------------------------------------------------
\section{Enumeration of CN- and PCN-elements}  \label{enumeration} 
%---------------------------------------------------------------------
\noindent 
In this section we explain the strategy underlying our Computational Result 3. 
It is based on the fundamental structure 
theory on completely normal elements from \cite{Ha-1997, Ha-2001a}. Although most of the details may also be found in the more 
recent survey article \cite{Ha-2013}, it is necessary to summarize  
the basic facts which are crucial for our computational enumeration of CN- and PCN-elements.

%--5.1-----------------------------------------------------------------------------------------------------------------
\subsection{Generalized cyclotomic modules and their complete generators} 
%--5.1-----------------------------------------------------------------------------------------------------------------
Consider once more the extension $E/F$ of Galois fields, corresponding to the pair $(q,n)$, and let again $p$ be the characteristic of these fields, 
and $\sigma$ the Frobenius automorphism of $E/F$. 
For a divisor $m$ of $n'$, let $\Phi_m(x)$ denote the $m$-th cyclotomic polynomial. 
A \textbf{generalized cyclotomic polynomial} (for $E/F$) has the form $\Phi_k(x^t)$, where 
$kt\mid n$ and $k$ is not divisible by $p$. Without loss of generality, one can additionally impose that $\gcd(k,t)=1$. 
%, but it 
%is convenient to allow a bit more flexibility. 
Since $\Phi_k(x^t)$ divides $x^n-1$,  
the set 
\[ C_{k,t} := \{ v \in E: \Phi_k(\sigma^t)v=0\} \] 
is a $\sigma$-invariant $F$-subspace of $E$; it is called the \textbf{(generalized) cyclotomic module} (of $E/F$) corresponding to $(k,t)$. 
Its \textbf{module-character} is the number 
$kt/\rad(k)$, with $\rad$ as explaind at the beginning of Section 2. The significance of the module character 
relies on the fact that $C_{k,t}$ is an $\fqm$-vector space for all $m$ dividing $kt/\rad(k)$. 

An important feature of any cyclotomic module $C_{k,t}$ is that it admits a \textbf{complete generator}, that is 
an element $v\in C_{k,t}$ such that 
\begin{equation}\label{eqn:CN-order-criterion} 
  \Ord_{q^{d}}(v) = \Phi_{\rad(k)} \Big(x^{\frac{kt}{\rad(k)d}}\Big) {\mbox{ for every }} d\mid \tfrac{kt}{\rad(k)}. \end{equation}  
In other words, $v$ is an element which simultaneously generates $C_{k,t}$ with respect to \textit{all} its module-structures arising from the intermediate fields of $E/F$ which 
act on $C_{k,t}$. 

%--5.2-----------------------------------------------------------------------------------------------------------------
\subsection{The Complete Decomposition Theorem} 
%--5.2----------------------------------------------------------------------------------------------------------------- 
It is intuitive that any decomposition of the additive group of $E$ into 
a direct sum of cyclotomic modules induces an additive decomposition of any completely normal element of $E/F$ 
into a sum of complete generators of the corresponding module components. The converse, however, is not true in general; it rather 
depends on the specific choice of the decomposition. 

In order to make this more precise, let us fix a cyclotomic module of $E/F$, 
say $C_{\ell,s}$ -- in the special case where $(\ell,s)=(1,n)$ this is just the extension field $E$ itself. 
Then a set $I$ of pairs $(k,t)$ is said to \textit{induce} a \textbf{cyclotomic decomposition} 
for $(\ell,s)$, respectively for $\Phi_\ell(x^s)$ and for $C_{\ell,s}$, 
provided that $\Phi_{k_1}(x^{t_1})$ and $\Phi_{k_2}(x^{t_2})$ are relatively prime for any two distinct pairs $(k_1,t_1)$ and $(k_2,t_2)$ of $I$, 
and $\prod_{(k,t) \in I} \Phi_k(x^t) = \Phi_{\ell}(x^s)$. 
According to this, we have 
\[ C_{\ell,s} = \bigoplus_{(k,t) \in I} C_{k,t}, \] 
and any $v\in C_{\ell,s}$ accordingly is additively decomposed as $v=\sum_{(k,t)\in I} v_{(k,t)}$. 
Moreover, if $v$ is a complete generator of $C_{\ell,s}$, then 
every $v_{(k,t)}$ is a complete  generator of $C_{k,t}$. 

Next, $I$ is said to induce an \textbf{agreeable decomposition} provided that, conversely, \textit{any} collection 
$(u_{(k,t)}: (k,t) \in I)$, with $u_{(k,t)}$ being a complete generator of $C_{k,t}$ for all $(k,t)$, 
gives that $\sum_{(k,t)\in I} u_{(k,t)}$ is a complete generator of $C_{\ell,s}$. 
We are now able to formulate the following fundamental result from \cite{Ha-1997, Ha-2001a}: 

\begin{quote} 
\textit{Complete Decomposition Theorem.} 
Consider a generalized cyclotomic module $C_{k,t}$, as part of a Galois field extension $E/F$ with characteristic $p$. 
Let $r$ be a prime divisor of $t$ and write $t=r^as$, with $s$ not dividsible by $r$. 
Assume that $r\not=p$ and that $r$ does not divide $k$. 
Then 
\[ I_r:= \left\{ \left(k,\tfrac{t}{r}\right),  \left(kr^a,\tfrac{t}{r^a}\right) \right\} \] 
induces a cyclotomic decomposition of $C_{k,t}$. 
Moreover, the following two statements are equivalent:  
\begin{enumerate} 
\item[(i)]  $I_r$ induces an agreeable decomposition of $C_{k,t}$;  
\item[(ii)]  $\ord_{\rad(kt')}(q)$ is not divisible by $r^a$.  
 \end{enumerate} 
\end{quote} 

\noindent 
The Complete Decomposition Theorem (DCT for short) is always applicable to $(1,n)$ with $r$ being the largest prime divisor of 
$n'$. Furthermore, it usually  may be applied iteratively several times. For instance, when $n=r^m$ is a prime power 
(with $r\not=p$), then the canoncial decomposition 
$x^n-1 = \prod_{i=0}^m \Phi_{r^i}(x)$ is agreeable. It is also important to note that the module character is \textit{reduced} by an application of DCT, 
namely from $kt/\rad(k)$ of the initial cyclotomic module to $\frac{1}{r} \cdot kt/\rad(k)$ of any of its two parts. 

%--5.3-----------------------------------------------------------------------------------------------------------------
\subsection{Enumerating CN-elements} 
%--5.3----------------------------------------------------------------------------------------------------------------- 
Since the process of a successive refinement of an aggreeable decomposition relying on DCT  
is confluent by \cite{Ha-2001a}, every cyclotomic module admits a \textbf{finest agreeable decomposition}. 

Throughout, we let $I_{q,n}^*$ denote the index set of the finest agreeable decomposition of the field extension $\fqn$ over $\fq$, and for 
every pair $(k,t)\in I_{q,n}^*$, we define $\phi_q^c[k,t]$ to be the total number of all 
complete generators of the cyclotomic module $C_{k,t}$ over $\fq$. 
As an immediate consequence, we have 
\begin{equation}\label{eqn:CN-number}  
CN_n(q) = \prod_{(k,t)\in I_{q,n}^*} \phi_q^c[k,t] . 
\end{equation} 
Let us have a look at a concrete situation. 

\begin{quote} 
\textit{Example} 5.3.1. 
When $q=3$ and $n=20$, then 
$\{ (1,1), (2,1),(4,1), (5,4)\}$ induces the finest agreeable decomposition of $\mathbb{F}_{3^{20}}$ over ${\mathbb F}_3$, namely 
\[ \mathbb{F}_{3^{20}} = {\mathbb F}_3 \oplus C_{2,1} \oplus C_{4,1} \oplus C_{5,4}, \] 
corresponding to $x^{20}-1 = \Phi_1(x) \Phi_2(x) \Phi_4(x)\Phi_5(x^4)$. 
The numbers of complete generators for these cyclotomic modules 
are as follows: 
\[ \phi^c_3[1,1] = 2 = \phi^c_3[2,1] \ {\mbox{ and }} \ \phi^c_3[4,1]=8 \ {\mbox{ and }} \ \phi^c_3[5,4]=37\, 015\, 040. \] 
With quation \eqref{eqn:CN-number} we achieve $CN_3(20)=1 \, 184 \, 481\, 280$.  %\qed  
\end{quote} 

\noindent In comparison to Morgan and
Mullen \cite{Mo-Mu-1996}, the use of DCT enables us to widen the range for enumerations of CN- and PCN-elements enormously. 
We shall outline the general approach, first, before subsequently emphazising several special aspects. 
%%the calculation 
%%Figure \ref{fig:alg:enumeration} gives an outline of the algorithmic
%%approach of  
%%We are going to explain the relevant steps now: 

\begin{enumerate} 
\item[(1)] Given a pair $(q,n)$, we start by determining {\em some} completely normal element $w$ for $\fqn$ over $\fq$ 
as described in Section 4 , in particular in Subsection 4.4.  
\item[(2)] We further determine the (index set $I_{q,n}^*$ of the) finest agreeable decomposition of $\fqn$ over $\fq$. 
\item[(3)] For every $(k,t)\in I_{q,n}^*$, let $\Gamma_{(k,t)}(x):=(x^n-1)/\Phi_k(x^t)$ and 
\[ u_{(k,t)}:=\Gamma_{(k,t)}(\sigma)w. \]   
The oberservation that $\Gamma_{(k,t)}(x)$ is equal to $(y^{\rad(k)}-1)/\Phi_{\rad(k)}(y)$, where $y=x^{kt/\rad(k)}$,  
yields that $u_{(k,t)}$ is a complete generator for the cyclotomic component $C_{k,t}$ of $\fqn$. 
\item[(4)] For every $(k,t)\in I_{q,n}^*$ we determine the number $\phi_q^c[k,t]$ of all 
complete generators of $C_{k,t}$ over $\fq$ from $u_{(k,t)}$. This is explained in detail below.  
\item[(5)] After that, we obtain $CN_n(q)$ from \eqref{eqn:CN-number} . 
\end{enumerate} 
Suppose, we are given some cyclotomic module $C_{k,t}$, where $(k,t)\in I_{q,n}^*$, and let us consider its (first) 
complete generator $u:=u_{(k,t)}$ defined in Step (3) above. 
Then, 
\[ C_{k,t} = \{ h(\sigma)u: h(x) \in \fq[x], \deg(h)< \varphi(k) t\}. \]  
Moreover, if $h(x) \in \fq[x]$ with $\deg(h)< \varphi(k)t$, then 
$h(\sigma)u$ has $q$-order equal to $\Phi_k(x^t)$ if and only if $h(x)$ and $\Phi_k(x^t)$ are relatively prime, which means that 
$h(x)$ corresponds to a unit in the residue ring $\fq[x]/(\Phi_k(x^t))$. 
We iterate through the set of all these $h(x)$ and at each time we 
check, whether $v:=h(\sigma)u$ satisfies the condition \eqref{eqn:CN-order-criterion}, where $d=1$ is already covered by the choice of $h(x)$ resulting in $v$. 
Let, for short, $\kappa:=kt/\rad(k)$ denote the module character of $C_{k,t}$. 
Whether $v$ has the correct $q^d$-order for all the remaining divisors $d$ of $\kappa$ can be 
performed in the same way as the complete normality test in Subsection 4.2: 

\begin{quote} 
for every $d\mid \kappa$ with $d\not=1$, consider the distinct monic irreducible divisors $g_1(x)$, ..., $g_s(x)$ of 
$\Phi_{\rad(k)}(x^{\kappa/d})$ over $\fqd$, and let $G_i(x):= \Phi_{\rad(k)}(x^{\kappa/d})/g_i(x)$ be their corresponding cofactors 
(for $i=1,...,s$); then $v$ has $q^d$-order equal to $\Phi_{\rad(k)}(x^{\kappa/d})$ if and only if 
$G_i(\sigma^d)v\not=0$ for every $i=1,...,s$. 
\end{quote} 

\begin{quote} 
\textit{Remark} 5.3.2. 
Advantage can be drawn from Proposition 3.2.4 as follows: let 
$D_{(k,t)}$ be the set of all divisors of $kt/\nu(k)$. Then $v$ is already a complete generator of $C_{k,t}$,  
when 
\begin{equation}\label{eqn:relaxed-order-conditon} 
 \Ord_{q^{d}}(v) = \Phi_{\rad(k)} \Big(x^{\frac{kt}{\rad(k)d}}\Big) \ {\mbox{ for every }} \ d \in {\mathcal D}^* \cap D_{(k,t)}, \end{equation}  
where ${\mathcal D}^*$ is the $(q,n)$-essential set associated to the CN-digraph for $(q,n)$. 
\end{quote} 

\noindent 
We emphasize that the concept of regularity can be generalized to 
cyclotomic modules, as well: $C_{k,t}$ is \textbf{regular} provided that $\ord_{\rad(kt')}(q)$ and $kt$ are relatively
prime. In that case, an element is a complete generator of $C_{k,t}$ over $\fq$ if 
it already has the correct $q^d$-order for at most \textit{two} specific members $d$ of $D_{(k,t)}$. 
For details and a summary we refer to  \cite[Section 20]{Ha-1997} and \cite[Section 5.6.4]{Ha-2013}. 

%--5.4-----------------------------------------------------------------------------------------------------------------
\subsection{Enumerating PCN-elements} 
%--5.4----------------------------------------------------------------------------------------------------------------- 
The derivation of the total number $PCN_n(q)$ of all primitive completely normal elements for some pair $(q,n)$ 
requires to additively recombine \textit{every} completely normal element from its cyclotomic components corresponding to $I_{q,n}^*$ 
and perform the primitivity test as described in Subsection 4.3. For the instance $(q,n)=(3,20)$, for example, 
we get 
\[ PCN_3(20)=423\, 266\, 160. \] 
Recall from Table \ref{tab:enumeration-DH-SH-1} and Table \ref{tab:enumeration-DH-SH-2} that our range comprises degrees $n$ which are less than $32$. 
Similar to the example where $q=3$ and $n=20$, it is therefore quite typical that 
$I_{q,n}^*$ produces one \textit{big} component, while all other components are \textit{small}. 
In the example just mentioned, $(1,1)$, $(2,1)$ and $(4,1)$ give the small components, while 
$(5,4)$ indicates the big one (as evident from the corresponding numbers of complete generators listed in Example 5.3.1).  
Based on this observation, during the process of enumerating $CN_n(q)$ and $PCN_n(q)$ for a given pair $(q,n)$, 
it turned out to be very fruitful 
to store \textit{all} complete generators for every small component in the memory of the computer, while 
dynamically generating the complete generators of the big component.

%---------------------------------------------------------------------
\section{Accessing the computational results}
\label{computational_results}
%---------------------------------------------------------------------
\noindent 
All software and resulting tables can be found under

\begin{itemize} 
\item 
\url{https://github.com/hackenbergstefan/Paper_PCN/}, and 
\item \url{https://github.com/hackenbergstefan/Masterarbeit/}. 
\end{itemize}

\subsection{Data for Computational Result 1} 
These are documented under 
\begin{quote} 
\url{https://github.com/hackenbergstefan/Paper_PCN/}. 
\end{quote}

\noindent 
In the folder \url{final} the tables 

\begin{itemize} 
\item \url{criterions_1_100.csv} ($n\leq 100$), \item \url{criterions_101_200.csv} 
($101\leq n \leq 200$), \item \url{criterions_201_202.csv} ($n\in \{201,202\}$). 
\end{itemize} 
can be found. According to what has been said in Subsection 3.3, each of these three tables consists of data of the following form:

\vspace{1em}

\begin{tabular}{|r|r|r|l|l|l|l|l|l|}
  \hline
  $p$ & $e$ & $n$ & C1 & C2 & C3 & C4 & C5 & C6 \\\hline\hline
  $2$ & $2$ & $10$ & False & False & False & False & False & $x^{20} + x^{19} + x^4 + x^3 + 1$\\\hline
  3 & 2 & 10 & False & False & False & False & True & \\\hline
  89 & 1 & 100 & True & True & True & & &\\\hline
\end{tabular}

\vspace{1em}

\noindent 
That is, these tables contain a line for each triple $(p, e, n)$ with $q := p^e < n'$ and $n \leq 202$, where the pair $(p^e, n)$ is not regular.
C1 to C5 represent the criterions given in Subsection \ref{sufficient-result-1}. 
An explicit PCN-polynomial is provided in column C6 if all other criterions fail.

\subsection{Data for Computational Result 2} These are also documented under 
\begin{quote} 
\url{https://github.com/hackenbergstefan/Paper_PCN/} 
\end{quote}

\noindent 
and can be found in the folder \url{final/range}. The tables there have a naming of the form 

\begin{itemize} 
\item \texttt{pcns\_$p$.csv} (where $p< 10\,000$ is a prime number). 
\end{itemize} 
The following exemplary table is an excerpt of two files:

\vspace{1em}

\begin{tabular}{|r|r|l|p{7cm}|}
  \hline
  $p$ & $n$ & poly & factorization \\\hline\hline
  101 & 5 & $x^5 + x^4 + 2$ & $2^2 \cdot 5^3 \cdot 31 \cdot 491 \cdot 1381$ \\\hline
  233 & 33 & $x^{33} + x^{32} + 6$ & $2^3 \cdot 7 \cdot 23 \cdot 29 \cdot 7789 \cdot 3148333 \cdot 4494621011 \cdot 3891196548493 \cdot 4581484617271 \cdot 18075348903971940081205337161$ \\\hline
\end{tabular}

\vspace{1em}

\noindent 
The column ``poly'' gives a PCN-polynomial of $\F_{p^{n}}$ over $\F_p$, where $p^n < 10^{80}$.
The factorization of $p^{n} - 1$ is given in the column ``factorization''.

\subsection{Data for Computational Result 3}
The results of the enumerations of CN- and PCN-elements can be found in the repository of \cite{Ha-2015c}, that is  

\begin{quote} 
\url{https://github.com/hackenbergstefan/Masterarbeit/}. 
\end{quote}
The relevant folder is \url{Tables/Enumerations}. The tables support the naming patterns

\begin{itemize} 
\item 
\texttt{enumerationsPCN\_P\_$p$.csv} (where $p\leq 43$ is the characteristic of a finite field $\F_q$ for prime powers $q$ as in Table 2, see Section 1), 
\item 
\texttt{enumerationsPCN\_N\_$n$.csv} (where $n\in \{3,4,6\}$). 
\end{itemize} 
The first of these patterns concern the data in Table 2, while the second ones cover the data from Table 3 (see Section 1). 
The typical content of these tables is given as an excerpt as follows: 

\vspace{2ex} 

\begin{tabular}{|r|r|r|r|l|l|p{5cm}|}
  \hline
  $q$ & $p$ & $r$ & $n$ & CN & PCN & gens \\\hline\hline
  2 & 2 & 1 & 30 & 111132000 & 55308540 & (1 1 2)*: 2\newline (3 1 2): 12\newline (5 1 2): 240\newline (15 1 2): 57600 \\\hline
\end{tabular}
\vspace{1em}

\noindent 
The columns ``CN'', respectively ``PCN'' contain the number of CN-, respectively PCN-elements for $\fqn$ over $\fq$. 
Observe that, in contrast to the notation used in the present work, the notation $q = p^r$ (instead of $p^e$) 
in accordance with \cite{Ha-2015c} 
is used in these tables. 

The column ``gens'' contains the concrete numbers of complete generators for the particular cyclotomic modules
occuring in a finest agreeable decomposition of $\fqn$ over $\fq$ (see Subsections 5.2 and 5.3). 
For instance, $(k, t', \pi): N$ means that $\phi^c_q[k,t] = N$, where $t=t'\cdot \pi$, with $\pi$ being a power of $p$, while $\gcd(p,t')=1$ 
(see Subsection 5.3). 
Whenever a (generalized) cyclotomic module is regular, this has been indicated by $(\ .\ )^\ast$.

\subsection{Involved software}
We used \texttt{sage/python} to implement the theoretical results stated in this paper.
The files can be found in the folder \url{ff_pcn} in 

\begin{quote} 
\url{https://github.com/hackenbergstefan/Paper_PCN/}. 
\end{quote}
For the factorizations of $q^n - 1$ we made use of \texttt{yafu}\footnote{\texttt{yafu} is an acronym for 
\emph{yet another factoring utility}. 
Although the integer factoring procedures of \texttt{Sage} are fast, we used the
so called \texttt{yafu-setup-package} from 
\url{https://github.com/KingBowser/yafu-setup-package} which contains 
all sources and a top level Makefile for all needed utilities.}, 
which provides the most powerful modern algorithms to factor integers in a
completely automated way optimized for multithreaded processing. 

A small \texttt{readme} with installation and usage instructions is also placed in the repository.

\vspace{3ex} 
\noindent 
\textbf{Acknowledgements.} This research is based on the second author's Master thesis \cite{Ha-2015c}, 
written under the supervision of the first author.

We thank the Leibniz-Rechenzentrum of the Bavarian Academy of Sciences\footnote{https://www.lrz.de/} 
which empowered us by computational capabilities.

\vspace{3ex} 

\end{document}